\documentclass{article}
\usepackage{amsfonts, mathrsfs, amsmath, amsthm, amssymb}
\newtheorem{def1}{Definition}[section]
\newtheorem{pro}{Proposition}[section]
\newtheorem{rem}{Remark}[section]

\newtheorem{cor} {Corollary}[section]
\usepackage[all]{xy} 
\begin{document}
\newtheorem{thm}{Theorem}[section]
\title{${\bf {\it Q}}$-${\bf SOB}$ as an epireflective hull in ${\bf {\it Q}}$-${\bf TOP_0}$}
\author {Sheo Kumar Singh\thanks{sheomathbhu@gmail.com}\\ {\em Department of Mathematics},\\ {\em Banaras Hindu University},\\ {\em Varanasi-221 005, India}\\\\ Arun K. Srivastava\thanks{arunksrivastava@gmail.com} \\ {\em Department of Mathematics}\\ {\em and}\\ {\em Centre for Interdisciplinary Mathematical Sciences},\\ {\em Banaras Hindu University},\\ {\em Varanasi-221 005, India}} 
\date{}
\maketitle 
\abstract We show that the epireflective hull of the $Q$-Sierpinski space in the category ${\bf {\it Q}}$-${\bf TOP_{0}}$ of $Q$-$T_0$-topological spaces is the category ${\bf {\it Q}}$-${\bf SOB}$ of $Q$-sober topological spaces.\\


\section{Introduction} For a given (but fixed) variety $\textbf{A}$ of $\Omega$-algebras and a fixed member $Q$ of $\textbf{A}$, S.A. Solovyov \cite{Solo} introduced the notion of a $Q$-topological space (and $Q$-continuous maps between them), providing thereby the category ${\bf {\it Q}}$-${\bf TOP}$ \footnote{recently, a characterization of the category ${\bf {\it Q}}$-${\bf TOP}$ has been given in \cite{SKS}} of such spaces. He also introduced the notions of $Q$-$T_0$-topological spaces, $Q$-sober topological spaces and $Q$-Sierpinski space. If ${\bf {\it Q}}$-${\bf SOB}$ denotes the category of $Q$-sober topological spaces, then Solovyov also showed implicitely that ${\bf {\it Q}}$-${\bf SOB}$ is reflective in ${\bf {\it Q}}$-${\bf TOP}$ (cf. Lemma $19$ of \cite{Solo}).

 In this note, (motivated by results in \cite{Nel, SK}) we have shown that ${\bf {\it Q}}$-${\bf SOB}$ is the epireflective hull of the $Q$-Sierpinski space in the category ${\bf {\it Q}}$-${\bf TOP_0}$ of $Q$-$T_0$-topological spaces. 

\section{Preliminaries}
For all undefined category theoretic notions used in this paper, \cite{AHS} may be referred. All subcategories used here are assumed to be full.

We begin by recalling \footnote{Most of the definitions in the preliminaries are given in \cite{SKS} also, we recall these\\ here for the sake of completeness.} the notions of $\Omega$-algebras and their homomorphisms; for details, cf. \cite{mane}, \cite{Solo}.\\
\begin{def1} Let $\Omega = (n_\lambda)_{\lambda\in I}$ be a class of cardinal numbers. 
\begin{itemize}
\item An {\bf $\Omega$-algebra} is a pair $(A, (\omega _{\lambda}^{A})_{\lambda\in I})$ consisting of a set $A$ and a family of maps $\omega _{\lambda}^{A}: A^{n_\lambda}\rightarrow A$. $B\subseteq A$ is called a {\bf subalgebra} of $(A, (\omega _{\lambda}^{A})_{\lambda\in I})$ if $\omega _{\lambda}^{A}((b_i)_{i \in n_\lambda})\in B$,  for every $\lambda \in I$ and every $(b_i)_{i \in n_\lambda}\in B^{n_\lambda}$. Given $S\subseteq A$, $\langle S \rangle$ denotes the subalgebra of $(A, (\omega _{\lambda}^{A})_{\lambda\in I})$ `generated by $S$', i.e., $\langle S \rangle$ is the intersection of all subalgebras of $(A, (\omega _{\lambda}^{A})_{\lambda\in I})$ containing $S$. {\rm (}In fact, $\langle S \rangle= \{\omega _{\lambda}^{A}({\langle s_i\rangle}_{i \in n_\lambda})$ $|$ $s_i \in S$ and $\lambda \in I\}${\rm )}.   

\item Given $\Omega$-algebras $(A, (\omega _{\lambda}^{A})_{\lambda\in I})$ and $(B, (\omega _{\lambda}^{B})_{\lambda\in I})$, a map $f: A\rightarrow B$ is called an {\bf $\Omega$-algebra homomorphism} provided that for every $\lambda \in I$, the following diagram  
$$\xymatrix{A^{n_\lambda}\ar[d]_{\omega _{\lambda}^{A}}\ar[r]^{f^{n_\lambda}}&B^{n_\lambda}\ar[d]^{\omega _{\lambda}^{B}}\\
            A\ar[r]_{f}&B}$$ commutes.\\
 Let ${\bf Alg}(\Omega)$ denote the category of $\Omega$-algebras and $\Omega$-algebra homomorphisms {\rm (}this category has products{\rm )}.
            
\item A {\bf variety} of $\Omega$-algebras is a full subcategory of ${\bf Alg}(\Omega)$ which is closed under the formation of products, subalgebras, and homomorphic images.\\
  
  \quad  {\rm Throughout this paper, $\Omega = (n_\lambda)_{\lambda\in I}$ denotes a fixed class of cardinal numbers, $\textbf{A}$ denotes a fixed variety of $\Omega$-algebras and $Q$ denotes a fixed member of $\textbf{A}$}.\\
  
   \quad Each function $f: X \rightarrow Y$ between sets $X$ and $Y$ gives rise to two functions $f^{\leftarrow}: 2^Y \rightarrow 2^X$ and $f^{\rightarrow}: 2^X \rightarrow 2^Y$, given by $f^{\leftarrow}(B) = \{x\in X$ $|$ $f(x)\in B\}$ and $f^{\rightarrow}(A) = \{f(x)$ $|$ $x\in A\}$, and also a function $f_{Q}^{\leftarrow}: Q^Y \rightarrow Q^X$, given by $f_{Q}^{\leftarrow}(\alpha) = \alpha \circ f$.
\item Given a set $X$, a subset $\tau$ of $Q^X$ is called a {\bf $Q$-topology} on $X$ if $\tau$ is a subalgebra of $Q^X$, in which case, the pair $(X, \tau)$ is called a {\bf $Q$-topological space}. 

\item Given two $Q$-topological spaces $(X, \tau)$ and $(Y, \eta)$, a {\bf $Q$-continuous function} from $(X, \tau)$ to $(Y, \eta)$ is a function $f: X\rightarrow Y$ such that $f_{Q}^{\leftarrow}(\alpha) \in \tau$, for every $\alpha \in \eta$.

\item Given a $Q$-topological space $(X, \tau)$ and $Y\subseteq X$, $(i_{Q}^{\leftarrow})^{\rightarrow}(\tau)$ $(= \{p\circ i$ $|$ $p \in \tau\})$ is called the ${\bf {\it Q}}$-{\bf subspace topology} on $Y$, where $i: Y \rightarrow X$ is the inclusion map. We shall denote the $Q$-subspace topology on $Y$ as $\tau _{Y}$.  

\item A $Q$-topological space $(X, \tau)$ is called $Q$-${\bf T_0}$ if for every distinct $x, y \in X$, there exists $p\in \tau$ such that $p(x)\neq p(y)$.          
\end{itemize}
\end{def1}

The meanings of homeomorphisms, embeddings, and products etc. for $Q$-topological spaces, are on expected lines.

\quad  Let ${\bf {\it Q}}$-${\bf TOP}$ denote the category of all $Q$-topological spaces and $Q$-continuous maps between them.\\

\begin{rem}
In {\rm \cite{Solo}}, it has been noted that, ${\bf {\it Q}}$-${\bf TOP}$, like ${\bf TOP}$, has products. One can go on further and verify that ${\bf {\it Q}}$-${\bf TOP}$ is initially complete; in fact ${\bf {\it Q}}$-${\bf TOP}$ turns out to be a topological category over {\rm \textbf{SET}}. As a consequence of the above, ${\bf {\it Q}}$-${\bf TOP}$ is complete; in particular, it has equalizers which are constructed, at the set-theoretical level, in the same way as in {\rm \textbf{SET}}.
\end{rem}
Let $\mathscr {C}$ be a category, $\mathscr {H} \subseteq mor\mathscr{C}$, $epi \mathscr{C}$ be the class of all $\mathscr C$-epimorphisms and $\mathscr R$ be a subcategory of $\mathscr {C}$. 
\begin{def1}{\rm \cite{AHS, Giu}}
$\mathscr R$ is said to be {\rm epireflective} in $\mathscr C$ if for each $\mathscr C$-object $X$, there exists an epimorphism $r_{X}:X\rightarrow RX$, with $RX\in ob\mathscr {R}$, such that for each $\mathscr C$-morphism $f:X\rightarrow Y$, with $Y\in ob\mathscr {R}$, there exists a unique $\mathscr {R}$-morphism $f^{*}:RX\rightarrow Y$, such that $f^{*}\circ r_{X}=f$. If moreover, each $r_{X}\in \mathscr {H}$ and $f^{*}$ is a $\mathscr C$-isomorphism, whenever $f\in epi\mathscr{C}\bigcap \mathscr H$, then $\mathscr R$ is said to be an $\mathscr {H}$-{\rm firm epireflective subcategory} of $\mathscr C$ {\rm (}or that the epireflectivity of $\mathscr R$ in  $\mathscr C$ is $\mathscr {H}$-firm{\rm )}.
\end{def1}
\section{The $Q$-sober space}
 Consider the identity function $id: Q\rightarrow Q$ and let $\nu = \langle id\rangle$ be the subalgebra of $Q^{Q}$, generated by $id$. 
\begin{def1} {\rm \cite{Solo}}
The $Q$-topological space $(Q, \nu)$ is called the $Q$-Sierpinski space.
\end{def1}
We shall denote the $Q$-Sierpinski space $(Q, \nu)$ as $Q_S$.\\\\
The next result is from \cite{SKS} (which is also the same as Lemma $57$ in \cite{Solo}).
\begin{thm}
For every $Q$-topological space $(X, \tau)$, $p\in \tau$ if and only if \\$p: (X, \tau)\rightarrow Q_{S}$ is $Q$-continuous.
\end{thm}
 For every $A \in ob\textbf{A}$, let $ptA= hom_{\textbf{A}}(A, Q)$. Define a map $\phi: A\rightarrow Q^{ptA}$ as $\phi(a)(p)= p(a)$, $\forall a\in A$ and $\forall p\in ptA$. Then $\phi$ turns out to be an $\Omega$-algebra homomorphism (cf. \cite{Solo}). Hence $\phi(A)$ is a subalgebra of $Q^{ptA}$, whereby $\phi(A)$ is a $Q$-topology on $ptA$.

 For each $Q$-topological space $(X, \tau)$, let $X\xrightarrow{\eta_{X}} pt\tau$ be the function defined by $\eta_{X}(x)(p)= p(x)$, $\forall x\in X$ and $\forall p\in \tau$.
\begin{def1} {\rm \cite{Solo}}
A $Q$-topological space $(X, \tau)$ is called $Q$-{\rm sober} if $\eta_{X}$ is bijective.
\end{def1} 
Let ${\bf {\it Q}}$-${\bf SOB}$ denote the subcategory of ${\bf {\it Q}}$-${\bf TOP}$ whose objects are $Q$-sober topological spaces.

The following fact is easily verified.
\begin{pro}
The $Q$-Sierpinski space $Q_{S}$ is $Q$-sober.
\end{pro}
\begin{pro} {\rm \cite{Solo}}
For every $A\in ob \textbf{A}$, $(ptA, \phi(A))$ is $Q$-sober.
\end{pro}


\begin{pro} {\rm \cite{Solo}}
Let $(X, \tau)\in ob{\bf {\it Q}}$-${\bf TOP}$. Then
\begin{enumerate}
\item the $Q$-topological space $(pt\tau, \phi(\tau))$ is $Q$-sober,
\item $\eta_{X}:(X, \tau)\rightarrow (pt\tau, \phi(\tau))$ is $Q$-continuous,
\item $(X, \tau)$ is $Q$-$T_{0}$ if and only if $\eta_{X}$ is injective,
\item $(X, \tau)$ is $Q$-sober if and only if $\eta_{X}:(X, \tau)\rightarrow (pt\tau, \phi(\tau))$ is $Q$-homeomorphism.
\end{enumerate}
\end{pro}
\begin{pro}
${\bf {\it Q}}$-${\bf SOB}$ is reflective in ${\bf {\it Q}}$-${\bf TOP}$.
\end{pro}
\textbf{Proof}: It follows from Lemma $19$ of \cite{Solo}.
\begin{pro}
If $(X, \tau)\in ob{\bf {\it Q}}$-${\bf TOP_0}$, then $\eta_{X}:(X, \tau)\rightarrow (pt\tau, \phi(\tau))$ is a ${\bf {\it Q}}$-${\bf TOP_0}$-embedding.
\end{pro}
\textbf{Proof}: Clearly, $\eta_{X}:(X, \tau)\rightarrow (pt\tau, \phi(\tau))$ is injective. Let $f:(X, \tau)\rightarrow (\eta_{X}(X), \phi(\tau)_{\eta_{X}(X)})$ be the `corestriction' of $\eta_{X}$ onto $\eta_{X}(X)$. It is enough to show that $f^{-1}: (\eta_{X}(X), \phi(\tau)_{\eta_{X}(X)})\rightarrow (X, \tau)$ is $Q$-continuous, i.e., to show that $p\circ f^{-1}\in \phi(\tau)_{\eta_{X}(X)}$, $\forall p\in \tau$, where $\phi(\tau)_{\eta_{X}(X)}$ is the $Q$-subspace topology on $\eta_{X}(X)$. Note that $\phi(\tau)_{\eta_{X}(X)} = \{\phi(q)\circ i$ $|$ $q\in \tau \}$, where $i: \eta_{X}(X)\rightarrow pt \tau$ be the inclusion map. For every given $p\in \tau$, as $(p\circ f^{-1})(\eta_{X}(x)) = p(f^{-1}(\eta_{X}(x))) = p(x)$ and also, $(\phi(p)\circ i)(\eta_{X}(x)) = \phi(p)(\eta_{X}(x)) = \eta_{X}(x)(p) = p(x)$, $\forall x\in X$, so, $p\circ f^{-1} = \phi(p)\circ i$, $\forall p\in \tau$. Hence $f^{-1}$ is $Q$-continuous.

\subsection{Another description of $Q$-sobriety}
This section is motivated by \cite{Pump} (Section 6; Prop. 28(c), page 107) wherein, sobriety in ${\bf TOP}$ was shown to have some link with an adjoint situation between ${\bf TOP}$ and $\textbf{SET}^{op}$, arising out of a use of the two-point Sierpinski topological space. We show here that an analogous link exists for $Q$-sobriety also.

Let $G: \textbf{SET}^{op}\rightarrow {\bf {\it Q}}$-${\bf TOP}$ and $F: {\bf {\it Q}}$-${\bf TOP}\rightarrow \textbf{SET}^{op}$ be the functors, described as follows:

    $G$ sends an object $X$ to $Q_{S}^X$ (the $X$-fold product of $Q_{S}$) and a morphism  $f:X\rightarrow Y$ to $G(f):Q_{S}^X\rightarrow Q_{S}^Y$, given by $G(f)(g) = g\circ f$, while $F$ sends an object $X=(X, \tau)$ to the set $C(X, Q_S)$ of all $Q$-continuous functions from $(X, \tau)$ to $Q_{S}$ (which is equal to $\tau$) and a morphism $f:X\rightarrow Y$ to $F(f):C(X, Q_S)\rightarrow C(Y, Q_S)$, given by $F(f)(\alpha)=\alpha\circ f$.

It can be easily verified that $G$ is right adjoint to $F$ and that the unit $\psi:Id_{Q\textbf{-TOP}}\rightarrow GF$ of this adjunction is given as follows:
 
for every $X=(X, \tau)\in ob{\bf {\it Q}}$-${\bf TOP}$, $\psi_X:X\rightarrow GFX (=Q_S^{\tau})$, is defined as $\psi_{X}(x)(p)= p(x), \forall x\in X, \forall p\in \tau$.

 Let $T= GF$ and let $e:EX\rightarrow TX$ be the equalizer (in ${\bf {\it Q}}$-${\bf TOP}$) of $T\psi_X$ and $\psi_{TX}$, $e$  being the inclusion map. Then $EX=\{f\in TX\mid T\psi_X(f)=\psi_{TX}(f)\}$. As $\psi$ is a natural transformation, we have $T\psi_X\circ\psi_X =\psi_{TX}\circ\psi_X$, so there exists a unique morphism $k_X: X\rightarrow EX$ in ${\bf {\it Q}}$-${\bf TOP}$ such that $\psi_X=e\circ k_X$ (see the following diagram). 

 $$\xymatrix
 {X\ar[dr]^{{\rm\psi_X}}\ar[d]_{{\rm k_X}}\\EX\ar[r]_e&TX\ar @< 2pt >[r]^{T\psi_X}\ar @< -2pt >[r]_{\psi_{TX}}&T^2X}$$
   
 \begin{pro}
 Let $(X,\tau)\in$ ob${\bf {\it Q}}$-${\bf TOP}$. Then $(X,\tau)$ is $Q$-sober iff $k_X: X\rightarrow EX$ is a $Q$-homeomorphism.
 \end{pro}
\noindent
  \textbf{Proof:} In view of Prop. $3.3(4)$, it will suffice to show that $(i)$ the $Q$-topological space $EX$ is the same as $(pt\tau, \phi (\tau))$ and $(ii)$ $k_X= \eta_X$. 

For $(X, \tau) \in$ ob${\bf {\it Q}}$-${\bf TOP}$, note that $T\psi_X, \psi_{TX}: Q_S^\tau\rightarrow Q_S^{C(Q_S^\tau, Q_S)}$ are given by $(T\psi_X)(f)(\alpha)=f(\alpha\circ \psi_X)$ and $(\psi_{TX})(f)(\alpha)=\alpha(f), \forall f\in Q_S^\tau$, $\forall\alpha\in C(Q_S^\tau, Q_S)$. So, $EX=\{f\in TX\mid f(\alpha\circ \psi_X)=\alpha(f), \forall\alpha\in C(Q_S^\tau, Q_S)\}$.  As $Q_S^\tau$ has the product $Q$-topology, so  for every $\alpha\in C(Q_S^\tau, Q_S)$, there is some $\lambda \in I$ such that $\alpha = {\omega_{\lambda}^Q}^{Q^{\tau}}({\langle \pi_{\alpha_i}\rangle}_{i\in n_{\lambda}})$, where $\pi_{\alpha_i}: {Q_{S}}^{\tau}\rightarrow Q_S$ is the ${\alpha_i}^{th}$ projection map, for $\alpha_i \in \tau$ and $i\in n_{\lambda}$. So, $\alpha \circ \psi_{X}= {\omega_{\lambda}^Q}^{Q^{\tau}}({\langle \pi_{\alpha_i}\rangle}_{i\in n_{\lambda}})\circ \psi_{X}= {\omega_{\lambda}^Q}^{X}({\langle \pi_{\alpha_i}\circ \psi_{X}\rangle}_{i\in n_{\lambda}})= {\omega_{\lambda}^Q}^{X}({\langle \alpha_i \rangle}_{i\in n_{\lambda}})$. Hence, $f(\alpha \circ \psi_{X})=f({\omega_{\lambda}^Q}^{X}({\langle \alpha_i \rangle}_{i\in n_{\lambda}}))$. Also, $\alpha (f)= ({\omega_{\lambda}^Q}^{Q^{\tau}}({\langle \pi_{\alpha_i}\rangle}_{i\in n_{\lambda}}))(f)= {\omega_{\lambda}^Q}({\langle \pi_{\alpha_i}(f)\rangle}_{i\in n_{\lambda}})= {\omega_{\lambda}^Q}({\langle f(\alpha_i)\rangle}_{i\in n_{\lambda}})$. Consequently, $EX= \{f\in TX$ $|$ $f({\omega_{\lambda}^Q}^{X}({\langle \alpha_i \rangle}_{i\in n_{\lambda}}))= {\omega_{\lambda}^Q}({\langle f(\alpha_i)\rangle}_{i\in n_{\lambda}}), \forall \lambda \in I\}$, i.e., $EX= \{f$ $|$ $f:\tau \rightarrow Q$ is an $\Omega$-algebra homomorphism$\} = pt\tau$.

Note also that $EX$ is the $Q$-subspace of $Q_{S}^{\tau} (=TX)$. It can be easily verified that the $Q$-subspace topology on $EX$ is the same as the $Q$-topology $\phi(\tau)$ on $pt\tau$. Thus the $Q$-topological spaces $EX$ and $(pt\tau, \phi (\tau))$ are the same. This establishes $(i)$.

From the definition of $\psi_X$, it is clear that $\eta_X$ is the `corestriction' of $\psi_X$ to $pt\tau$. Also, from the diagram above, $\psi_X=e\circ k_X$.  Hence, $\forall x\in X$, $\eta_X(x)=\psi_X(x)=e(k_X(x))=k_X(x)$. Thus $\eta_X= k_X$, which establishes $(ii)$.  $\Box$   
    
\subsection{${\bf {\it Q}}$-${\bf SOB}$ as the epireflective hull of $Q_{S}$}

For $(X, \tau)\in ob {\bf {\it Q}}$-${\bf TOP}$ and $M\subseteq X$, put $[M]=\bigcap\{Eq(f,g)$ $|$ $f,g \in \tau$ and $f|_{M}=g|_{M}\}$, where $Eq(f,g)=\{x\in X$ $|$ $f(x)=g(x)\}$. It turns out that $[[M]]=[M]$. Also, if $[M]=M$, then we say that $M$ is $[\quad]$-closed.\\ 

For showing that ${\bf {\it Q}}$-${\bf SOB}$ is the epireflective hull of $Q_{S}$ in ${\bf {\it Q}}$-${\bf TOP_0}$, we shall need to identify $(i)$ the epimorphisms in ${\bf {\it Q}}$-${\bf TOP_0}$ and $(ii)$ the extremal subobjects in ${\bf {\it Q}}$-${\bf TOP_0}$.
 
\begin{pro}
A morphism $e:(X,\tau)\rightarrow (Y,\delta)$ in ${\bf {\it Q}}$-${\bf TOP_0}$ is an epimorphism if and only if $e_{Q}^{\leftarrow}$ is injective.
\end{pro}
\textbf{Proof}: Suppose $e$ is an epimorphism and for $q_{1},q_{2}\in \delta$, $e_{Q}^{\leftarrow}(q_1)=e_{Q}^{\leftarrow}(q_2)$. Then $q_{1}\circ e=q_{2}\circ e$, implying that $q_1=q_2$.

Conversely, suppose the given condition is satisfied. Now, consider any distinct pair $f,g:(Y, \delta)\rightarrow (Z, \sigma)$ of morphisms in ${\bf {\it Q}}$-${\bf TOP_0}$. Then for some $y\in Y$, $f(y)\neq g(y)$. Since $Z$ is $Q$-$T_{0}$, $\exists$ $p\in \sigma$ such that $p(f(y))\neq p(g(y))$, i.e., $f_{Q}^{\leftarrow}(p)\neq g_{Q}^{\leftarrow}(p)$. This gives $e_{Q}^{\leftarrow}(f_{Q}^{\leftarrow}(p))\neq e_{Q}^{\leftarrow}(g_{Q}^{\leftarrow}(p))$ i.e., $p\circ f\circ e\neq p\circ g\circ e$, implying that $f\circ e\neq g\circ e$. Thus $e$ is an epimorphism. $\Box$
\begin{pro}
A morphism $f:(X, \tau)\rightarrow (Y, \delta)$ in ${\bf {\it Q}}$-${\bf TOP_0}$ is an epimorphism if and only if $[f(X)]=Y$.
\end{pro}
\textbf{Proof}:\footnote{This result has been proved in \cite{Cast} (Theorem $1.11$) in a more general set-up. The\\ proof being given here is somewhat more direct.} First, let $f:(X, \tau)\rightarrow (Y, \delta)$ be an epimorphism in ${\bf {\it Q}}$-${\bf TOP_0}$. Let $[f(X)]\neq Y$. Then $\exists$ $ $ $y\in Y$ such that $y\notin [f(X)]$, and so $\exists$ morphisms $g,h:(Y, \delta)\rightarrow Q_{S}$ in ${\bf {\it Q}}$-${\bf TOP_0}$ with $g|_{f(X)}= h|_{f(X)}$ and $g(y)\neq h(y)$. Since $g|_{f(X)}= h|_{f(X)}$, $g\circ f=h\circ f$, which is a contradiction. Thus $[f(X)]= Y$.

Conversely, let $[f(X)]= Y$. Consider any two  morphisms $g, h:(Y, \delta)\rightarrow (Z, \sigma)$ in ${\bf {\it Q}}$-${\bf TOP_0}$ such that $g\circ f=h\circ f$. If possible, let $g\neq h$. Then $\exists$ $ $ $y\in Y$ such that $g(y)\neq h(y)$. Since, $g\circ f=h\circ f$, $g|_{f(X)}= h|_{f(X)}$. But then $y\notin [f(X)]$, a contradiction. Thus $f$ is an epimorphism. $\Box$\\

We say that an embedding $e:(X, \tau)\rightarrow (Y, \delta)$ in ${\bf {\it Q}}$-${\bf TOP_0}$ is $[\quad]$-closed if $[e(X)]=e(X)$.
\begin{pro}
The extremal monomorphisms in ${\bf {\it Q}}$-${\bf TOP_0}$ are precisely the $[\quad]$-closed embeddings $($in ${\bf {\it Q}}$-${\bf TOP_0})$.
\end{pro}
\textbf{Proof}:  Let $m:(X, \tau)\rightarrow (Y, \delta)$ be an extremal monomorphism in ${\bf {\it Q}}$-${\bf TOP_0}$ and $Z=[m(X)]$ (with the $Q$-subspace topology $\delta_{Z}$). Define a map $e:(X, \tau)\rightarrow (Z, \delta_{Z})$ as $e(x)=m(x)$, $\forall x\in X$. Then $e$ is an epimorphism in ${\bf {\it Q}}$-${\bf TOP_0}$ and $m=i\circ e$, where $i:(Z, \delta_{Z})\rightarrow (Y, \delta)$ is the inclusion map. But then $e$ is a $Q$-homeomorphism. Thus $m$ is a $[\quad]$-closed embedding.

Conversely, let $m:(X, \tau)\rightarrow (Y, \delta)$ be a $[\quad]$-closed embedding. Let the elements of the set $\{(f, g) \in \tau \times \tau$ $|$ $f|_{m(X)}=g|_{m(X)}\}$ be indexed by an index set $J$. Then $[m(X)]=\bigcap\{Eq(f_j, g_j)$ $|$ $f_j, g_j \in \tau$ and $j \in J\}$. For every $j\in J$, let $\pi_j: Q_{S}^{J} \rightarrow Q_{S}$ be the $j^{th}$ projection map. Then by the property of the product, there exists unique $Q$-continuous maps $f^*, g^*: (X, \tau)\rightarrow  Q_{S}^{J}$ such that $\pi_j \circ f^* = f_j$ and $\pi_j \circ g^* = g_j$, $\forall j\in J$. Now it can be easily verified that $[m(X)] = Eq (f^*, g^*)$, whereby $m(X)$ (in fact, the inclusion map from $(m(X), \delta_{m(X)})$ to $(Y, \delta)$) is an equalizer in ${\bf {\it Q}}$-${\bf TOP_0}$ (this also follows from Proposition $1.6$ of \cite{Cast}, which, however, is stated in a more general set-up). But as equalizers are extremal monomorphisms, $m$ is an extremal monomorphism. $\Box$

\begin{cor}
The extremal subobjects in ${\bf {\it Q}}$-${\bf TOP_0}$ are precisely the $[\quad]$-closed subspaces of $Q$-$T_{0}$-topological spaces.
\end{cor}

The next result is analogous to the corresponding results in \cite{Nel} and \cite{SK}.
\begin{thm}$(i)$ ${\bf {\it Q}}$-${\bf SOB}$ is epireflective in ${\bf {\it Q}}$-${\bf TOP_0}$ and $(ii)$ this epireflectivity is $\mathscr {H}$-firm, where $\mathscr {H}$ is the class of all ${\bf {\it Q}}$-${\bf TOP_0}$-embeddings.
\end{thm}
\textbf{Proof}:
 $(i)$ Let $(X,\tau)\in ob{\bf {\it Q}}$-${\bf TOP_0}$. We show that $\eta_{X}: (X,\tau)\rightarrow (pt\tau, \phi(\tau))$ is the desired epireflection of $(X, \tau)$ in ${\bf {\it Q}}$-${\bf SOB}$. We use Proposition $3.6$ to show first that $\eta_{X}$ is an epimorphism. Let ${\eta_{X}}_{Q}^{\leftarrow}(\phi(p_{1}))={\eta_{X}}_{Q}^{\leftarrow}(\phi(p_{2}))$, where $p_{1}, p_{2}\in \tau$. Then ${\eta_{X}}_{Q}^{\leftarrow}(\phi(p_{1}))(x)={\eta_{X}}_{Q}^{\leftarrow}(\phi(p_{2}))(x)$, $\forall x\in X$, implying that $\phi(p_{1})(\eta_{X}(x))=\phi(p_{2})(\eta_{X}(x))$, i.e., $\eta_{X}(x)(p_{1})=\eta_{X}(x)(p_{2})$, which gives $p_{1}(x)=p_{2}(x)$. Hence, $p_{1}=p_{2}$, whereby $\phi(p_{1})=\phi(p_{2})$. 

Let $(Y, \delta)\in ob{\bf {\it Q}}$-${\bf SOB}$ and $f:(X, \tau)\rightarrow(Y, \delta)$ be $Q$-continuous. We need to find a ${\bf {\it Q}}$-${\bf TOP}$-morphism $f^{*}:(pt\tau, \phi(\tau))\rightarrow (Y,\delta)$ such that $f^{*}\circ \eta_{X}=f$. For any $\alpha\in pt\tau$, define $\alpha':\delta \rightarrow Q$ by $\alpha'(q)=\alpha(q\circ f)$. It can be verified that $\alpha'$ is an $\Omega$-algebra homomorphism, i.e., $\alpha'\in pt\delta$. Since $(Y, \delta)$ is $Q$-sober, there is a unique $y\in Y$ with $\eta_{Y}(y)=\alpha'$. Put $f^{*}(\alpha)= y$. This gives us a map $f^{*}:pt\tau \rightarrow Y$. Now, given $q \in \delta$ and $\alpha\in pt\tau$, ${f^{*}}_{Q}^{\leftarrow}(q)(\alpha)=(q\circ f^{*})(\alpha)=q(f^{*}(\alpha))=q(y)=\eta_{Y}(y)(q)=\alpha'(q)=\alpha(q\circ f)=\alpha(f_{Q}^{\leftarrow}(q))=\phi(f_{Q}^{\leftarrow}(q))(\alpha)$, whereby ${f^{*}_{Q}}^{\leftarrow}(q)=\phi(f_{Q}^{\leftarrow}(q))$. Hence ${f^{*}}_{Q}^{\leftarrow}(q)\in \phi(\tau)$, showing the $Q$-continuity of $f^{*}$. Now, $\forall q \in \delta$ and $\forall x\in X$, $\eta_{Y}(f(x))(q)=q(f(x))=(q\circ f)(x)=\eta_{X}(x)(q\circ f)=(\eta_{X}(x))'(q)$. Hence, $\eta_{Y}(f(x))=(\eta_{X}(x))'$, $\forall x\in X$. So $f^{*}(\eta_{X}(x))=f(x)$, $\forall x\in X$. Hence, $f^{*}\circ \eta_{X}=f$. Finally, as $\eta_{X}$ is an epimorphism, $f^{*}$ is unique. This proves $(i)$.

We now prove $(ii)$. Let $(X, \tau)\in ob{\bf {\it Q}}$-${\bf TOP_0}$. Then the epireflection $\eta_{X}: (X, \tau)\rightarrow (pt\tau, \phi(\tau))$ is injective and it is an embedding in ${\bf {\it Q}}$-${\bf TOP_0}$. Let $(Y, \delta)\in ob{\bf {\it Q}}$-${\bf SOB}$, $f:(X, \tau)\rightarrow(Y, \delta)$ be an epimorphic-embedding in ${\bf {\it Q}}$-${\bf TOP_0}$ and $f^{*}:(pt\tau, \phi(\tau))\rightarrow (Y, \delta)$ the unique ${\bf {\it Q}}$-${\bf TOP}$-morphism such that $f^{*}\circ \eta_{X}=f$. Let $\tilde{f}: X \rightarrow f(X)$ be the `corestriction' of $f$ to $f(X)$. Then ${\tilde{f}}^{-1}: (f(X), \delta_{f(X)})\rightarrow (X, \tau)$ is clearly $Q$-continuous. So, $\forall p\in \tau$ $ $ $\exists$ $ $ $p_{f}\in \delta$ such that $(\tilde{f}^{-1})_{Q}^{\leftarrow}(p)= p_{f}\circ i$, where $i: f(X)\rightarrow Y$ is the inclusion map. Hence $f_{Q}^{\leftarrow}(p_f)=p$. As $f$ is an epimorphism, this $p_f$ is unique such that $f_{Q}^{\leftarrow}(p_f)=p$. Now, define $g:(Y, \delta)\rightarrow (pt\tau, \phi(\tau))$ by $g(y)(p)= p_{f}(y)$, $\forall y\in Y$ and $\forall p\in \tau$. It can be  easily verified that $g(y): \tau \rightarrow Q$ is an $\Omega$-algebra homomorphism. Now, $\forall p \in \tau$ and $\forall y\in Y$, $g_{Q}^{\leftarrow}(\phi(p))(y)=(\phi(p)\circ g)(y)=\phi(p)(g(y))=g(y)(p)=p_{f}(y)$. So, $g_{Q}^{\leftarrow}(\phi(p))=p_{f}$. Thus $g_{Q}^{\leftarrow}(\phi(p))\in \delta$, showing that $g$ is $Q$-continuous. Since, $\forall q \in \delta$, $f_{Q}^{\leftarrow}(q)\in \tau$ and so $f_{Q}^{\leftarrow}((f_{Q}^{\leftarrow}(q))_{f})=f_{Q}^{\leftarrow}(q)$, whereby $(f_{Q}^{\leftarrow}(q))_{f}=q$. Hence, $\forall q \in \delta$ and $\forall y\in Y$, $(f_{Q}^{\leftarrow}(q))_{f}(y)=q(y)$, implying that $g(y)(f_{Q}^{\leftarrow}(q))=\eta_{Y}(y)(q)$, i.e., $g(y)(q\circ f)=\eta_{Y}(y)(q)$. So, $f^{*}(g(y))=y$, $\forall y\in Y$. Thus $f^{*}\circ g=id_{Y}$.

Next, let $\alpha \in pt\tau$ and $f^{*}(\alpha)= y$. Then $\eta_{Y}(y)(q)=\alpha(q\circ f)$, $\forall q \in \delta$. For $p\in \tau $, $g(y)(p)=p_{f}(y)=\eta_{Y}(y)(p_{f})=\alpha (p_{f}\circ f)=\alpha (f_{Q}^{\leftarrow}(p_{f}))=\alpha(p)$ implying that $g(y)=\alpha$. Hence $g\circ f^{*}=id_{pt\tau}$. Thus $f^{*}$ is a ${\bf {\it Q}}$-${\bf TOP}$-isomorphism. $\Box$

Using Theorem $1$ of \cite{Mar}, we get the following corollary:
\begin{cor}
${\bf {\it Q}}$-${\bf SOB}$ is closed under forming products and extremal subobjects in ${\bf {\it Q}}$-${\bf TOP_0}$.
\end{cor}

\begin{pro}
If $(X, \tau)\in ob {\bf {\it Q}}$-${\bf TOP_0}$, then $(pt\tau, \phi(\tau))$ is a $[\quad]$-closed subspace of $Q_{S}^{\tau}$.
\end{pro}
\textbf{Proof}: From \cite{Solo} (Theorem $58$), it follows that the map $e:X\rightarrow Q_{S}^{\tau}$ defined by $e(x)(\mu)=\mu(x)$, $\forall x\in X$, $\forall \mu\in \tau$, is a ${\bf {\it Q}}$-${\bf TOP_0}$-embedding. Hence $f: X \rightarrow [e(X)]$, the `corestriction' of $e$ to $[e(X)]$, is an epimorphic-embedding in ${\bf {\it Q}}$-${\bf TOP_0}$ (by Proposition $3.7$). $[e(X)]$, being a $[\quad]$-closed subspace of $Q_{S}^{\tau}$ (which is $Q$-sober), is, therefore, an extremal subobject of $Q_{S}^{\tau}$. Hence $[e(X)]$ is $Q$-sober. Theorem $3.2(ii)$, now provides a $Q$-sober isomorphism $f^{*}:pt\tau \rightarrow [e(X)]$ (such that $f^{*}\circ \eta_{X}= f$). Hence $(pt\tau, \phi(\tau))$ is a $[\quad]$-closed subspace of $Q_{S}^{\tau}$. $\Box$
\begin{pro}
$(X, \tau)\in ob{\bf {\it Q}}$-${\bf SOB}$ if and only if it is $Q$-homeomorphic to a $[\quad]$-closed subspace of $Q_{S}^{\tau}$.
\end{pro}
\textbf{Proof}: Let $(X, \tau)\in ob{\bf {\it Q}}$-${\bf SOB}$. Then, via $\eta_{X}$, $(X, \tau)$ is $Q$-homeomorphic to $(pt\tau, \phi(\tau))$, which is a $[\quad]$-closed subspace of $Q_{S}^{\tau}$. The converse follows from Proposition $3.1$, Corollary $3.1$ and Corollary $3.2 $. $\Box$\\

Using Theorem $2$ of \cite{Mar}, together with Corollary $3.1$ and Proposition $3.10$ above, we now obtain the following result:
\begin{thm}
${\bf {\it Q}}$-${\bf SOB}$ is the epireflective hull of $Q_{S}$ in ${\bf {\it Q}}$-${\bf TOP_0}$.
\end{thm}

\textit{Acknowledgement}: The first author would like to thank the Council of Scientific and Industrial Research, New Delhi, India, for financial support through its Senior Research Fellowship.

\end{document}